\documentclass[11pt,a4paper]{article}
\usepackage{amsfonts,amsmath,amssymb,amstext,amsbsy,amsopn,amsthm}
\usepackage{url}
\usepackage{graphicx}

\newcommand{\beq}{\begin{eqnarray}}
\newcommand{\eeq}{\end{eqnarray}}
\newcommand{\beqn}{\begin{eqnarray*}}
\newcommand{\eeqn}{\end{eqnarray*}}

\newcommand{\itg}{\int \limits}

\def\C{{\mathbb C}}
\newcommand{\Ree}{\mbox{\rm Re}\:}

\newcommand{\cA}{\mathcal A}

\newcommand{\cD}{\mathcal D}

\newcommand{\cK}{\mathcal K}

\newcommand{\cM}{\mathcal M}

\newcommand{\cO}{\mathcal O}
\newcommand{\cP}{\mathcal P}
\newcommand{\cQ}{\mathcal Q}

\newcommand{\cS}{\mathcal S}

\newcommand{\cF}{\mathcal F}
\newcommand{\cE}{\mathcal E}

\newcommand{\bdx}{\mathbf{x}}
\newcommand{\bdy}{\mathbf{y}}
\newcommand{\bdm}{\mathbf{m}}

\newcommand{\bdk}{\mathbf{k}}
\newcommand{\bdxi}{\boldsymbol{\xi}}

\DeclareMathOperator{\erfc}{erfc}

\DeclareMathOperator{\e}{e}

\def\R{{\mathbb R}}
\def\Z{{\mathbb Z}}

\newcommand\ds{\displaystyle}
\def\l{\lambda}

\def\bdz{{\bf z}}

\newtheorem{thm}{Theorem}[section]
\newtheorem{lem}{Lemma}[section]
\newtheorem{rem}{Remark}[section]

\title{Computation of volume potentials over bounded domains\\
via approximate approximations}

\author{{  F. Lanzara$^{\mbox{\tiny 1}}$ , V. Maz'ya$^{\mbox{\tiny 2}}$ ,
G. Schmidt$^{\mbox{\tiny 3}}$}
}
\date{}

\begin{document}
\maketitle

\hspace*{1cm}
\parbox{10cm}{\begin{flushleft}
{\footnotesize\em
\begin{itemize}
\item[$^{\mbox{\tiny\rm 1}}$] Department of Mathematics, Sapienza University of Rome,\\
Piazzale Aldo Moro 2, 00185 Rome, Italy\\
\texttt{\rm lanzara\symbol{'100}mat.uniroma1.it}
\item[$^{\mbox{\tiny\rm 2}}$]Department of Mathematics, University of
Link\"oping, \\ 581 83 Link\"oping, Sweden;\\
Department of Mathematical Sciences, M\&O Building, University of
Liverpool, Liverpool L69 3BX, UK;\\
\texttt{\rm vlmaz\symbol{'100}mai.liu.se }
\item[$^{\mbox{\tiny\rm 3}}$]Weierstrass Institute for Applied Analysis and
Stochastics, \\  Mohrenstr. 39,
10117 Berlin, Germany \\
\texttt{\rm schmidt\symbol{'100}wias-berlin.de}
\end{itemize}
}
\end{flushleft}}

{\bf  Abstract.}  We obtain cubature formulas of volume potentials over bounded domains combining the basis functions introduced in the theory of approximate approximations with  their integration  over the tangential-halfspace. Then the computation is reduced to the quadrature of one dimensional integrals over the halfline. We conclude the paper providing numerical tests which show that  these formulas  give very accurate approximations and confirm the predicted order of convergence.\section{Introduction}
\setcounter{equation}{0}

We consider  the volume potential of modified Helmholtz operators
in
$\R^n$
\[
\cA_n=-\Delta+\l^2
\]
with $\l\in\C$.
If $\l^2\neq 0$, then the fundamental solution $\kappa_\l(\bdx)$ is given by
\begin{align*}
\kappa_\l(\bdx)=
\frac{1}{(2 \pi)^{n/2}}
\Big( \frac{|\bdx|}{\lambda}\Big)^{1-n/2}K_{n/2-1} (\lambda|\bdx|) \, ,
\end{align*}
where $\l\in\C\setminus(-\infty,0]$ and $K_\nu$ is the modified Bessel
function of the second kind. The fundamental solutions of the Laplacian
($\l =0$) are well-known
\begin{align*}
\kappa_0(\bdx)
=\left\{
\begin{array}{ll}
\ds \frac{1}{2\pi} \log \frac{1}{|\bdx|} \, , &n=2\, ,\\
\ds \frac{\Gamma(\frac{n}{2}-1)}
    {4 \pi^{n/2}} \frac{1}{|\bdx|^{n-2}}\, ,
&n\geq 3\, .
\end{array}
\right.
\end{align*}

For $f\in C^1(\Omega)$,  $\Omega\subset\R^n$, the volume potential
\begin{equation}\label{vol}
u(\bdx)=\cK_\l f (\bdx)=\int_{\Omega}\kappa_{\l}(\bdx-\bdy)f(\bdy)d\bdy
\end{equation}
  provides a solution of
\[
\cA_n u=\left\{ \begin{array}{cc} f(\bdx)& \bdx\in\Omega\\ \\ 0& otherwise
\end{array}\right.
\]
We study cubature formulas for	the volume potential \eqref{vol}  using
the concept {\it approximate approximations}  (see  \cite{MSbook}).
The cubature of  volume potentials over the full space and over
high-dimensional halfspaces  has already been studied in \cite{LMS2} and
\cite{LMS}, respectively. In \cite{MSW} cubature formulas based on approximate approximations for the single layer harmonic  potential  were considered.

Assume $f\in C^N(\Omega)$. We extend $f$ outside $\Omega$  with preserved
smoothness and we denote by  $\widetilde{f} \in C^N_0(\R^n)$ the continuation of
$f$.  Assume that there exists $C>0$ such that
\[
||\widetilde{f}||_{W^N_\infty(\R^n)}\leq C\, ||f||_{W^N_\infty(\Omega)}.
\]
We introduce a uniform grid $\{h\bdm\}$ with step $h$.
A cubature formula for \eqref{vol} can be obtained if we replace $f$ by the
approximate quasi-interpolant
\begin{equation}\label{quasiint}
\cM_{h,\cD}\widetilde{f}(\bdx)=\cD^{-n/2}
\sum_{\bdm\in\Z^n}\widetilde{f}(h\bdm)\eta\left(\frac{\bdx-h\bdm}
{h\sqrt{\cD}}\right)
\end{equation}
where  $\eta\in\cS(\R^n)$ and satisfies the moment conditions of order $N$
\begin{equation}\label{moment}
\int_{\R^n}\eta(\bdx)\,\bdx^{\alpha}d\bdx=\delta_{0,\alpha},\quad 0\leq
|\alpha|<N.
\end{equation}
The quasi-interpolant \eqref{quasiint} approximates $f$ in
$\Omega$.
It is known (\cite{MSbook}) that
\begin{equation*}
|{f}(\bdx)-\cM_{h,\cD}\widetilde{f}(\bdx)| \le
c (\sqrt{\cD}h)^N \|\nabla_N f\|_{L^\infty}
+ \sum_{k=0}^{N-1} \varepsilon_k (\sqrt{\cD}h)^k \big|\nabla_k f(\bdx)\big|
\end{equation*}
with
\begin{align*}
\varepsilon_k \le\sum_{\bdm \in {\Z}^n \setminus \{0\}}
  \big| \nabla_k \cF\eta (\sqrt{\cD} \bdm)  \big| \, ; \lim_{\cD\to \infty}
\sum_{\bdm \in {\Z}^n \setminus \{0\}}
  \big| \nabla_k \cF\eta (\sqrt{\cD} \bdm)=0.
\end{align*}
Since $\eta$ is a smooth and rapidly decaying function,
for any error $\epsilon >0$ one can fix $r>0$ and the parameter $\cD>0$ such
that the quasi-interpolant with nodes in a neighborhood of $\Omega$
\begin{equation}\label{quasiintr} 
\cM^r_{h,\cD}\widetilde{f}(\bdx)=
\cD^{-n/2}\sum_{\substack{d(h\bdm,\Omega)\leq
r\,h\sqrt{\cD}}}
\widetilde{f}(h\bdm)\eta\left(\frac{\bdx-h\bdm}
{h\sqrt{\cD}}\right)
\end{equation}
approximates $f$ with
\begin{equation}\label{estimate}
|{f}(\bdx)-\cM^r_{h,\cD}\widetilde{f}(\bdx)| = \cO((\sqrt{\cD}h)^N+ \epsilon)
\|f\|_{W^N_\infty}
\end{equation}
for all $x \in \Omega$.

We use the quasi-interpolant \eqref{quasiintr} of the density to obtain a cubature formula of the volume potential \eqref{vol}
\begin{equation}\label{cub}
\cK_{\l,h}\widetilde{f}(\bdx)=\cK_\l(\cM^r_{h,\cD}\widetilde{f}) (\bdx)
=
\cD^{-n/2}\!\!\!\! \sum_{\substack{\bdm\in\Z^n\\d(h\bdm,\Omega)\leq
r\,h\sqrt{\cD}}}
\widetilde{f}(h\bdm) \int_\Omega
\kappa_{\l}(\bdx-\bdy)\,
\eta\left(\!\frac{\bdy-h\bdm}{h\sqrt{\cD}}\right)d\bdy.
\end{equation}
Since $\cK_\l$ is a bounded mapping between suitable function spaces, the
differences
$\cK_{\l,h}\widetilde{f}(\bdx)-\cK_{\l}f(\bdx)$ behave like estimate
\eqref{estimate}. We use the
radial generating functions
\begin{equation}\label{defeta}
\eta_{2M}(\bdx)=\pi^{-n/2} L_{M-1}^{(n/2)}(|\bdx|^2) \e^{-|\bdx|^2}\, ,
\end{equation}
 where  $L_{k}^{(\gamma)}$ are the generalized Laguerre polynomials
     \begin{equation} \label{defLag}
L_{k}^{(\gamma)}(y)=\frac{\e^{\,y} y^{-\gamma}}{k!} \, \Big(
{\frac{d}{dy}}\Big)^{k} \!
\left(\e^{\,-y} y^{k+\gamma}\right), \quad \gamma > -{\rm 1} \, ,
     \end{equation}
which satisfy the moment condition
\eqref{moment} with $N=2M$. 
Formula \eqref{cub}  would give a cubature of \eqref{vol}, if the integrals  \[
  \int_\Omega  \kappa_{\l}(\bdx-\bdy)\,
\eta\left(\!\frac{\bdy-h\bdm}{h\sqrt{\cD}}\right)d\bdy
\]
could be computed efficiently 
for nodes with $d(h\bdm,\Omega)\leq
r\,h\sqrt{\cD}$.

 For smooth  domains, we propose to replace the integrals in \eqref{cub} by integrals over the {\it tangential-halfspace}  at a point of $\partial \Omega$ with minimal distance to $h \bdm$.  It is proved that  these formulas approximate \eqref{vol}  with the order $\cO((h\sqrt{\cD})^2)$.  We conclude the paper providing numerical tests which show that  these formulas  give very accurate approximations and the order of convergence cannot be improved.

\section{Cubature based on  \eqref{quasiint}}\label{sectwo}
\setcounter{equation}{0}

In this section we study the approximation of  the integral \eqref{vol}
over a bounded region $\Omega\subset\R^n$ with smooth boundary
by the sum \eqref{cub} for appropriately chosen $r>0$.

Denote by $P_{h\bdm}$ a point of $\partial \Omega$ with minimal distance to $h \bdm$  and by $n_{P_{h\bdm}}$ the normal at $P_{h\bdm}$ directed towards the interior of $\Omega$.  Let $T_{h\bdm}$ be the halfspace  bounded by the
tangential  plane  at  the point $P_{h\bdm}$ such that  
the inner normal  at $P_{h\bdm}$ coincides with  $n_{P_{h\bdm}}$. 

We define the following cubature formula for the volume potential
 \eqref{vol}
\begin{equation} \label{cub1}
\widetilde{\cK}_{\l,h}\widetilde{f}(\bdx)=\cD^{-n/2}
\sum_{d(h\bdm,\Omega)\leq r\,h\sqrt{\cD}}
\widetilde{f}(h\bdm) \int_{T_{h\bdm}}	\kappa_{\l}(\bdx-\bdy)\,
\eta\left(\!\frac{\bdy-h\bdm}{h\sqrt{\cD}}\right)d\bdy.
\end{equation}
To derive explicit expressions for the integrals
in \eqref{cub1},  as basis functions we use \eqref{defeta},
  satisfying the moment condition
\eqref{moment} with $N=2M$ (cf. \cite{MSbook}). 
We have
\begin{lem}\label{halfspace}
Let
$\Ree   \l^2\ge 0$, $n \ge 1$ or $\Ree   \l^2= 0$, $n \ge 3$.
The solution of the equation
\begin{equation*}
(-\Delta+\l^2) u=\left\{
\begin{tabular}{ccc}
$  \eta_{2M}(\bdx)$,&& $x_n\geq a$,\\
\\
0, &&$x_n<a$
\end{tabular}
  \right.
\end{equation*}
is given by the one-dimensional integral
\begin{align} \label{int_tang}
\frac{1}{8\pi^{n/2}} \int_0^\infty \e^{-  \l^2 t/4}\e^{-|\bdx|^2/(1+  t)}
\Big(\erfc\big(F(  t,x_n,a)\big) \cP_M(\bdx,  t)
+\frac{ \e^{-F^2(  t,x_n,a)}}{\sqrt{\pi} }
\cQ_M(\bdx,  t,a)\Big)\,dt
\end{align}
where
\begin{align*}
\cP_M(\bdx,t)&=\sum_{k=0}^{M-1}\frac{1}{(1+t)^{k+n/2}} 
L_k^{(n/2-1)}\big(\frac{|\bdx|^2}{1+t}
\big) \, ,\\
\cQ_M(\bdx,t,a)
&=\frac{2}{({1+t})^{(n-1)/2}}
\sum_{k=0}^{M-1} \sum_{l=0}^k \frac{(-1)^{k-l}}{(k-l)!
4^{k-l}}L_l^{((n-3)/2)}
\Big(\frac{|\bdx'|^2}{1+t}\Big)
\sum_{j=1}^{2(k-l)}\frac{(-1)^j}{t^{j/2}} \\
&\times\bigg( \!
\Big(\hskip-3pt\begin{array}{c}2(k-l)\\j\end{array}\hskip-3pt\Big)
H_{2(k-l)-j}\Big(\frac{x_n}{\sqrt{1+t}}\Big)
\frac{H_{j-1}(F(t,x_n,a))}{(1+t)^{k+1/2}}
-H_{j-1}\Big(\frac{a-x_n}{\sqrt{t}}\Big)
\frac{H_{2(k-l)-j}(a)}{(1+t)^l} \!\bigg)
\end{align*}
with the function 
\begin{align} \label{defF}
F(t, x, a)  =
  \sqrt{\frac{1+t}{t}} \Big( a -
  \frac{ x}{1+t} \Big).
\end{align}
$H_k$  denote the Hermite polynomials
\begin{equation}\label{defHerm}
H_k(x)=(-1)^k \e^{x^2}
\frac{d^k}{dx^k}
\e^{-x^2} \, .
\end{equation}
\end{lem}
\begin{proof}
We consider
the heat equation in $\R^n$
\begin{equation*}
\partial_t z-   \Delta z=0 \, , \quad
z(\bdx,0)=\left\{
\begin{tabular}{ccc}
$    \eta_{2M}(\bdx)$,&& $x_n\geq a$,\\
\\
0, &&$x_n<a$,
\end{tabular}
  \right.
\end{equation*}
whose solution
is given by the Poisson  integral
\begin{align*}
z(\bdx,t) =
\frac{1}{{(4\pi   t)^{n/2}}}
\itg_{y_n>a}\e^{-|\bdx-\bdy|^2/(4  t)}\eta_{2M}(\bdy)\,d\bdy.
\end{align*}
The representation  (\cite[p.55]{MSbook})
\begin{equation}\label{gen}
\eta_{2M}(\bdx)=\pi^{-n/2} \sum_{j=0}^{M-1} \frac{(-1)^j}{j!4^j}\Delta^j
\e^{-{|\bdx|^2}},
\end{equation}
 shows that
\[
\eta_{2M}(\bdx',x_n)=\pi^{-n/2}\sum_{j=0}^{M-1}\frac{(-1)^j}{j!4^j}
\sum_{l=0}^j \Big(\begin{array}{c}j\\l\end{array}\Big)
\Delta_{\bdx'}^l\e^{-|\bdx'|^2}\frac{d^{2(j-l)}}{d x_n^{2(j-l)}}\e^{-x_n^2}.
\]
Hence
\begin{align*}
z(\bdx,t) =\frac{1}{(4   t)^{n/2}\pi^{n}}
\sum_{j=0}^{M-1}\frac{(-1)^j}{j!4^j}
\sum_{l=0}^j \Big(\begin{array}{c}j\\l\end{array}\Big)
  \varphi_{j-l}(x_n,4  t,a)\> \int_{\R^{n-1}}\e^{-|\bdx'-\bdy'|^2/(4  t)}
\Delta_{\bdy'}^l \e^{-|\bdy'|^2}d\bdy'
\end{align*}
with 
\begin{align} \label{defvarphi}
\varphi_k(x, t,p) = \itg_{p}^{\infty} {\e}^{-(x- y)^2/t}
  \, \frac{d^{2k}}{dy^{2k}}  \e^{-y^2} \, dy.
   \end{align}
From
\begin{align*}
\int_{\R^{n-1}}\e^{-|\bdx'-\bdy'|^2/(4  t)} \Delta_{\bdy'}^l
\e^{-|\bdy'|^2}d\bdy'
&= \Delta_{\bdx'}^l\int_{\R^{n-1}}\e^{-|\bdx'-\bdy'|^2/(4  t)}
\e^{-|\bdy'|^2}d\bdy'
\\ &
=\Big(\frac{4 \pi   t}{1+4  t}\Big)^{(n-1)/2} \Delta^l_{\bdx'}
\e^{-|\bdx'|^2/(1+4  t)}
\end{align*}
and  the relation
\[
\Delta^j \e^{-|\bdx'|^2/(1+t)}=\frac{(-1)^j\,j
!4^j}{(1+t)^j}\e^{-|\bdx'|^2/(1+t)}L_j^{((n-3)/2)}
\Big(\frac{|\bdx'|^2}{1+t}\Big)
\]
(\cite[p.121]{MSbook}), we obtain
\begin{align*}
z(\bdx,t)& =\frac{1 }{(4   t)^{n/2}\pi^{n}}
\frac{\e^{-|\bdx'|^2/(1+4  t)}}{(1+4  t)^{(n-1)/2}}\sum_{j=0}^{M-1}\frac{(-1)^j
}{j!4^j} \\& \quad \times
\sum_{l=0}^j \Big(\hskip-3pt\begin{array}{c}j\\l\end{array}\hskip-3pt\Big)
  \frac{(-1)^l\,l !4^l}{(1+4  t)^l} \varphi_{j-l}(x_n,4  t,a)
L_l^{((n-3)/2)}
\Big(\frac{|\bdx'|^2}{1+4  t}\Big).
\end{align*}
From
\begin{align*}
\varphi_{0}(x,t,p)=\itg_{p}^{\infty} \e^{-(x-y)^2/t} \e^{-y^2}    \,
dy=\frac{\sqrt{\pi}}{2}\sqrt{\frac{t}{1+t}}\e^{-x^2/(1+t)}
\erfc\big(F(t,x,p)\big) \, ,
\end{align*}
for $k \geq 1$, integration by parts   leads to
\begin{align*}
\itg_{p}^{\infty} \e^{-(x- y)^2/t}
\frac{d^{2k}}{dy^{2k}}\e^{-y^2}  \, dy_j
=\frac{\partial^{2k}}{\partial x^{2k}}\varphi_{0}(x,t,p)
-\sum_{\ell=0}^{2k-1}(-1)^\ell
\frac{\partial^{\ell}}{\partial y^{\ell}}
\e^{-(x - y)^2/t}
\frac{d^{2k-\ell-1}}{d y^{2k-\ell-1}}\e^{-y^{2}}
\Bigg|_{y=p} \, ,
\end{align*}
and the definition \eqref{defHerm}
gives
\begin{align*}
&\frac{d^{2k-\ell-1}}{d y^{2k-\ell-1}}\e^{-y^{2}}
=(-1)^{2k-\ell-1}
\e^{-y^{2}} H_{2k-\ell-1}(y) \, ,
\; \frac{\partial^{\ell}}{\partial y^{\ell}}
\e^{-(x - y)^2/t}
=\frac{(-1)^\ell\e^{-(x - y)^2/t}}{t^{\ell/2}} \,
H_\ell \Big( \frac{y-x}{\sqrt{t}}\Big) .
\end{align*}
In view of
\[
  \frac{d^{\ell}}{dx^{\ell}}\erfc(x)= \frac{2}{\sqrt{\pi}}(-1)^{\ell}
\e^{-x^2} H_{\ell-1}(x)
)
,\quad \ell \geq 1 \, ,
\]
one gets for $\ell<2k$
\begin{align*}
\frac{\partial^{2k-\ell}}{\partial x^{2k-\ell}} \erfc\big(F(t,x,p)\big)&=
\frac{(-1)^{2k-\ell}}{(t(1+t))^{k-\ell/2}}
\left[ \frac{d^{2k-\ell}}{dz^{2k-\ell}}\erfc(z)\right]_{z=F(t,x,p)}\\
&=\frac{2 \e^{-F^2(t,x,p)} }{\sqrt{\pi}(t(1+ t))^{k-\ell/2}} \,
H_{2k-\ell-1}(F(t,x,p)) \,.
\end{align*}
Therefore,  since
\[
\frac{d^\ell}{d x^\ell} \e^{-x^2/(1+t)}
=\frac{(-1)^\ell\e^{-x^2/(1+t)}}{(1+t)^{\ell/2}}
  H_\ell\Big(\frac{x}{\sqrt{1+t}}\Big)\, ,
\]
we obtain
\begin{align*}
\frac{\partial^{2k}}{\partial x^{2k}}& \, \varphi_{0}(x,t,p)
=  \frac{\sqrt{\pi t}}{2}\frac{\e^{-x^2/(1+t)}}{(1+ t)^{k+1/2}}
H_{2k}\Big(\frac{x}{\sqrt{1+t}}\Big)\erfc(F(t,x,p))\\
&\; -\frac{ \sqrt{t} \e^{-x^2/(1+t)} \e^{-F^2(t,x,p)}}{(1+ t)^{k+1/2}}
  \sum_{\ell=0}^{2k-1}
\Big(\hskip-4pt\begin{array}{c}2k\\\ell\end{array}\hskip-4pt\Big)
\frac{(-1)^{\ell}}{t^{k-\ell/2}}
H_\ell\Big(\frac{x}{\sqrt{1+t}}\Big)
H_{2k-\ell-1}(F(t,x,p)) \, .
\end{align*}
Thus simple transformations give
\begin{align}
& \varphi_{k}(x,t,p)= \e^{-x^2/(1+t)} \bigg(\erfc\big(F(t,x,p)\big)
H_{2k}\Big(\frac{x}{\sqrt{1+t}}\Big)\frac{\sqrt{\pi\, t}}{2(1+t)^{k+1/2}}
\label{varphik} \\
&+\e^{-F^2(t,x,p)}
  \sum_{\ell=1}^{2k}\frac{(-1)^{\ell}}{t^{(\ell-1)/2}}
\bigg( \! \Big(\hskip-3pt\begin{array}{c}2k\\\ell\end{array}\hskip-3pt\Big)
H_{2k-\ell}\Big(\frac{x}{\sqrt{1+t}}\Big)
\frac{H_{\ell-1}(F(t,x,p))}{(1+t)^{k+1/2}}
-H_{\ell-1}\Big(\frac{p-x}{\sqrt{t}}\Big)
H_{2k-\ell}(p) \!\bigg)\bigg).\nonumber
\end{align}

Using \eqref{varphik}, the relations
\[
H_{2j}(x)=(-1)^j \, 4^j\, j!\, L_j^{(-1/2)}(x^2)\quad \mbox{and} \quad
L_j^{(a+b+1)}(x+y)=\sum_{l=0}^j L_l^{(a)}(x)L_{j-l}^{(b)}(y)
\]
we find
\begin{align*}
&\sum_{j=0}^{M-1}\frac{(-1)^j}{j! 4^j} \frac{1}{(1+t)^{j+n/2}}\sum_{l=0}^j
\left(\!\!\!\begin{tabular}{c}$j$\\$l$\end{tabular}\!\!\!\right){(-1)^l\,l
!4^l}L_l^{((n-3)/2)}
\Big(\frac{|\bdx'|^2}{1+t}\Big)
H_{2(j-l)}\Big(\frac{x_n}{\sqrt{1+t}}\Big)
\\
&=\sum_{j=0}^{M-1} \frac{1}{(1+t)^{j+n/2}}\sum_{l=0}^j
L_l^{((n-3)/2)}\Big(\frac{|\bdx'|^2}{1+t}\Big)
L_{j-l}^{(-1/2)}\Big(\frac{x_n^2}{{1+t}}\Big)=\cP_M(\bdx,t)
\end{align*}
leading to
\begin{align*}
z(\bdx,t)=\frac{  \e^{-|\bdx|^2/(1+4  t)}}{2\pi^{n/2}} \Big(
\erfc\big(F(4  t,x_n,a)\big) {\cP_M(\bdx,4   t)}
+ \frac{\e^{-F^2(4  t,x_n,a)}}{\sqrt{\pi}}
\cQ_M(\bdx,4   t,a) \Big) .
\end{align*}
\end{proof}
In the particular case $n=2$
\begin{align*}
&\cP_1(\bdx,t)=\frac{1}{1+t} \, ; \quad
\cP_2(\bdx,t)=\frac{1}{1+t}\Big(
1+\frac{1}{1+t}-\frac{|\bdx|^2}{(1+t)^2}\Big) \, ; \\
&\cP_3(\bdx,t)=\cP_2(\bdx,t)+\frac{1}{1+t}\Big(\frac{1}{(1+t)^2}-2
\frac{|\bdx|^2}{(1+t)^3}+
\frac{|\bdx|^4}{2\,(1+t)^4}\Big)\\
& \qquad 
\qquad=\frac{1}{1+t}\Big(1+\frac{1}{1+t}-\frac{|\bdx|^2}{(1+t)^2}+\frac{1}{(1+t)
^2}-2 \frac{|\bdx|^2}{(1+t)^3}+
\frac{|\bdx|^4}{2\,(1+t)^4}\Big) \, ;  \\
&\cQ_1(\bdx,t,a)=0; \quad
\cQ_2(\bdx,t,a)=-  \frac{\sqrt{t}}{(t+1)^{3/2}}\big(
a+\frac{x_2}{1+t}\big);\\
&\cQ_3(\bdx,t,a)=
\frac{1}
     {4}
{\frac{\sqrt{t}}{(1 + t)^{3/2}}}\,
     \left( \frac{-2\,a\,t}{{\left( 1 + t \right) }} +
       \Big( a + \frac{{x_2}}{1 + t} \Big) \,
        \Big(
	 \frac{4\,{{|\bdx|}}^2-2\,x_2^2}
	  {{\left( 1 + t \right) }^2}-	 \frac{7}{{1 + t  }}+
{2\,a^2-5}\Big)  \right).
  	       \end{align*}

\begin{rem} \label{remark2}
For sufficiently large $|a|>r$ the integrands
in \eqref{int_tang} are approximated by
\begin{align*}
\left\{
\begin{array}{cc}
0 \, ,& a  \ge r \, ,\\
2 e^{-  \l^2
t/4}\e^{-|\bdx|^2/(1+ t)}
\cP_M(\bdx,t) , &a  \le -r  \, ,
\end{array} \right.
\end{align*}
with the error $\cO(\e^{-a^2})$, which is in accordance with
\begin{align*}
  \int_{\R^n} \kappa_{\l}(\bdx-\bdy) \eta_{2M}(\bdy) \, d \bdy
=\frac{1}{4\pi^{n/2}}\int_0^\infty  {\e^{-  \l^2
t/4}\e^{-|\bdx|^2/(1+  t)}}  \cP_M(\bdx,  t)\,dt .
\end{align*}
\end{rem}

\begin{thm}\label{mainTheorem} Suppose that the generating function $\eta\in\cS(\R^n)$
satisfies the moment condition \eqref{moment}
with $N\geq 2$.
Then the integral \eqref{vol}
is approximated by the sum
\begin{equation} \label{cub2}
\widetilde{\cK}_{\l,h}\widetilde{f}(\bdx)=\cD^{-n/2}
\sum_{d(h\bdm,\Omega)\leq r\,h\sqrt{\cD}}
\widetilde{f}(h\bdm) \int_{T_{h\bdm}}	\kappa_{\l}(\bdx-\bdy)\,
\eta\left(\!\frac{\bdy-h\bdm}{h\sqrt{\cD}}\right)d\bdy
\end{equation}
 with the error estimate
\begin{equation}\label{mainestimate}
|\widetilde{\cK}_{h,\l}\widetilde{f}
(\bdx)-{\cK}_{\l}f(\bdx)|=\cO(\epsilon+c_1\,
(h\sqrt{\cD})^2)
\end{equation}
provided $\partial \Omega$ has $C^2-$smoothness.
 The saturation term $\epsilon$ can be negligible small if $\cD$ is large enough.

\end{thm}
\begin{proof}
We  study the difference 
 $$\widetilde{\cK}_{\l,h}\widetilde{f}(\bdx)-\cK_{\l,h}\widetilde{f}(\bdx)
=\cD^{-n/2}
\sum_{d(h\bdm,\Omega)\leq r\,h\sqrt{\cD}}
\widetilde{f}(h\bdm) \int_{T_{h\bdm}}	\kappa_{\l}(\bdx-\bdy)\,
\eta\left(\!\frac{\bdy-h\bdm}{h\sqrt{\cD}}\right)d\bdy$$
$$-\cD^{-n/2}
\sum_{d(h\bdm,\Omega)\leq r\,h\sqrt{\cD}}
\widetilde{f}(h\bdm) \int_{\Omega}	\kappa_{\l}(\bdx-\bdy)\,
\eta\left(\!\frac{\bdy-h\bdm}{h\sqrt{\cD}}\right)d\bdy
 $$
as $h\to0$.

Since $\eta(\bdy)$ is supported in the small neighborhood  of the origin $|\bdy|\leq r$, if $h\bdm\in\Omega$ such that $d(h\bdm,\partial \Omega)>r h \sqrt{\cD}$ we replace each integral  
\[
 \int_{\Omega}	\kappa_{\l}(\bdx-\bdy)\,
\eta\left(\!\frac{\bdy-h\bdm}{h\sqrt{\cD}}\right)d\bdy
\]
by integrals over $T_{h\bdm}$ without loss accuracy.

Assume $\bdm\in\Z^n: d(h\bdm, \partial \Omega)\leq r h \sqrt{\cD}$.
 We choose a local coordinate system with the origin ${\bf 0}$ at the point $h\bdm$  and  the normal at the nearest point of $\partial \Omega$ from ${\bf 0}$ 
 coincides with the $x_n-$ axis $({\bf 0}',1)$, ${\bf 0'}\in\R^{n-1}$. 
In these coordinates system the halfspace  is defined  by ${T}={T}_{h\bdm}=\{\bdy=(\bdy',y_n): \bdy'\in\R^{n-1},y_n > \rho_{h\bdm}\}$ where
\[
\rho=\rho_{h\bdm}=\left\{
\begin{tabular}{cc}
$dist({\bf 0'},\partial \Omega)$& if ${\bf 0} \notin \Omega$\\
$-dist({\bf 0'},\partial \Omega)$& if ${\bf 0} \in \Omega$.
\end{tabular}
\right.
\] 
We have  $|\rho|< r h \sqrt{\cD}$.
We assume that in a neighborhood $U$ of the point $({\bf 0'},\rho)\in \partial \Omega$ the domain $\Omega$ is given by
\begin{equation}\label{phi}
y_n\geq \varphi(\bdy'),\quad  \varphi({\bf 0})=\rho,\> \nabla \varphi({\bf 0})=0.
\end{equation}
 
We introduce the change of variable $\bdz=\bdz(\bdy)$ defined as
\[\bdz'=\bdy',\>z_n=y_n-\varphi(\bdy')+\rho\]
 and  its inverse $\bdy=\bdy(\bdz)$
\[
\bdy'=\bdz',\> y_n=z_n +\varphi(\bdz')-\rho.
\]

Choose $h$ such that $B_{h r \sqrt{\cD}}\cap \Omega\subset U$ and  assume $r_0\geq r$ such that $\bdz(\Omega\cap B_{h r \sqrt{\cD}})=T\cap B_{h r_0 \sqrt{\cD}}$.
Then the difference
\[ \int_{{T}}	\kappa_{\l}((\bdx-h \bdm)-\bdy)\,
\eta\left(\!\frac{\bdy}{h\sqrt{\cD}}\right)d\bdy-\int_{\Omega}	\kappa_{\l}((\bdx-h \bdm)-\bdy)\,
\eta\left(\!\frac{\bdy}{h\sqrt{\cD}}\right)d\bdy
\]
takes the form
\begin{equation*}\label{diff}
\int_{T\cap B_{r_0 h \sqrt{\cD}}}	\left[\kappa_{\l}((\bdx-h \bdm)-\bdz)\,
\eta\left(\!\frac{\bdz}{h\sqrt{\cD}}\right)- 	\kappa_{\l}((\bdx-h \bdm)-\bdy(\bdz))\,
\eta\left(\!\frac{\bdy(\bdz)}{h\sqrt{\cD}}\right)\right]d\bdz
\end{equation*}
\begin{equation}\label{eq2}
=\int_{{T}\cap B_{r_0 h \sqrt{\cD}}(\bf{0})}\kappa_{\l}((\bdx-h \bdm)-\bdz))\,\left[\eta\left(\!\frac{\bdz}{h\sqrt{\cD}}\right)-
\eta\left(\!\frac{\bdy(\bdz)}{h\sqrt{\cD}}\right)\right]d\bdz
\end{equation}
\begin{equation}\label{eq1}
+\int_{{T}\cap B_{r_0 h \sqrt{\cD}}(\bf{0})}\left[\kappa_{\l}((\bdx-h \bdm)-\bdz)-\kappa_{\l}((\bdx-h \bdm)-\bdy(\bdz)\right]
\eta\left(\!\frac{\bdy(\bdz)}{h\sqrt{\cD}}\right)d\bdz.
\end{equation}
Consider the integral in \eqref{eq2}. In view of \eqref{phi} we can consider locally
 \[
 \varphi(\bdz')=\rho+\frac{1}{2}\,K\bdz'\cdot \bdz',\quad K=\{\partial_{ij}\varphi({\bf 0})\}_{i,j=1}^{n-1}.
 \]
Therefore we have
\begin{equation}\label{zy(z)}
|\bdz-\bdy(\bdz)|=|\varphi(\bdz')-\rho|\leq c_1\,(h\sqrt{\cD})^2,\>\> \bdz\in B_{r_0 h \sqrt{\cD}}(\bf{0})
\end{equation}
which leads to
\begin{equation}\label{eq4}
\left|\eta\left(\!\frac{\bdz}{h\sqrt{\cD}}\right)-
\eta\left(\!\frac{\bdy(\bdz)}{h\sqrt{\cD}}\right)\right|\leq
c_2\, h \,\sqrt{\cD}.
\end{equation}

Suppose  $|\bdx-h \bdm| \geq 2 r_0 h \sqrt{\cD}$. Since $|(\bdx-h \bdm)-\bdz|\geq |\bdx-h \bdm|/2$,  we obtain
\[
\int_{T\cap B_{r_0 h \sqrt{\cD}}(\bf{0})}\kappa_{\l}((\bdx-h \bdm)-\bdz)\,d\bdz \leq c \frac{(h\sqrt{\cD})^n}{|\bdx-h \bdm|^{n-2}}.
\]

If  $|\bdx-h \bdm| < 2 r_0 h \sqrt{\cD}$  then, for $a>r_0$,
\[
\int_{T\cap B_{r_0 h \sqrt{\cD}}(\bf{0})}\kappa_{\l}((\bdx-h \bdm)-\bdz)\,d\bdz \]
\[\leq c \int_{B_{a h \sqrt{\cD}}(\bdx-h\bdm) }
 \frac{d \bdz}{|(\bdx-h \bdm)-\bdz|^{n-2}} +c\,\int_{B_{r_0 h \sqrt{\cD}}({\bf{0}})\setminus B_{a h \sqrt{\cD}}(\bdx-h\bdm)}
 \frac{d \bdz}{|(\bdx-h \bdm)-\bdz|^{n-2}}.
\]
Obviously the first integral in the right-hand side is $\cO((h\sqrt{\cD})^2)$. To estimate the second integral we use the relation
\begin{equation}\label{maindis}
|\bdx-\bdz|\geq c(|\bdx'-\bdz'|+ h\sqrt{\cD}),\> \forall \bdz\in \R^n \setminus B_{ah\sqrt{\cD}}(\bdx),
\end{equation}
which implies (\cite[(4.2.13)]{MS})
\[
\int_{B_{r_0 h \sqrt{\cD}}({\bf{0}}) \setminus B_{ah\sqrt{\cD}}(\bdx)} \frac{d\bdz}{|\bdx-\bdz|^{n-2}}\]
\[\leq
c\int_{B_{r_0 h \sqrt{\cD}}({\bf{0}}) \setminus B_{ah\sqrt{\cD}}(\bdx)} \frac{d \bdz}{(|\bdx'-\bdz'|+ h\sqrt{\cD})^{n-2}}
\leq 
c (h\sqrt{\cD})^2. 
\]
Thus, for all $\bdx \in \R^n$,
\[
\int_{T\cap B_{r_0 h \sqrt{\cD}}(\bf{0})}\kappa_{\l}((\bdx-h \bdm)-\bdz)\,d\bdz \leq c \frac{(h\sqrt{\cD})^n}{(|\bdx-h \bdm|+h\sqrt{\cD})^{n-2}}.
\]

Therefore, if  $S_h$ denotes  the strip $S_h=\{\bdx\in\R^n: d(\bdx,\partial \Omega)\leq r h \sqrt{\cD}\}$, keeping in mind \eqref{eq4} we have
\[\sum_{h\bdm\in S_h}  |\widetilde{f}(h\bdm)| \int_{{T}\cap B_{r_0 h \sqrt{\cD}}}\kappa_{\l}((\bdx-h \bdm)-\bdz)\,
\left|\eta\left(\!\frac{\bdz}{h\sqrt{\cD}}\right)-
\eta\left(\!\frac{\bdy(\bdz)}{h\sqrt{\cD}}\right)\right|d\bdz\]
  \[
 \leq C\, ||{f}||_{\infty}(h\sqrt{\cD})^{n+1} \sum_{\substack{h\bdm\in S_h} }
  \frac{1}{(|\bdx-h \bdm|+h\sqrt{\cD})^{n-2}}\]
 \begin{equation}\label{estimateh2}
\leq  C\, ||{f}||_{\infty}(h \sqrt{\cD}) \int_{S_h}  
  \frac{d \bdy }{(|\bdx-\bdy|+h\sqrt{\cD})^{n-2}}\leq  
  C\,||f||_\infty ( h\, \sqrt{\cD})^2.
\end{equation}

\bigskip

It remains to consider the integral in \eqref{eq1}. 
We use the inequality  (see \cite[p.80]{MS})
\begin{equation}\label{dis100}
|\kappa_{\l}((\bdx-h \bdm)-\bdz)-\kappa_{\l}((\bdx-h \bdm)-\bdy(\bdz)|
\end{equation}
\[\leq\
c\, (h\sqrt{\cD})^2
\left[
\frac{1}{|(\bdx-h\bdm)-\bdz|^{n-1}}+\frac{1}{|(\bdx-h\bdm)-\bdy(\bdz)|^{n-1}}
\right]
\]
and obtain
\[
\int_{{T}\cap B_{r h \sqrt{\cD}}({\bf 0})}
\left|\kappa_{\l}((\bdx-h \bdm)-\bdz)-\kappa_{\l}((\bdx-h \bdm)-\bdy(\bdz)\right|
|\eta\left(\!\frac{\bdy(\bdz)}{h\sqrt{\cD}}\right)|\,d\bdz
\]
\[
\leq 
c_1\, (h\sqrt{\cD})^2
\int_{B_{r_0 h \sqrt{\cD}}({\bf 0})}
\frac{1}{|(\bdx-h\bdm)-\bdz|^{n-1}}d\bdz.
\]

If $|\bdx-h \bdm| \geq 2 r_0 h \sqrt{\cD}$ then
\[
\int_{ B_{r_0 h \sqrt{\cD}}(\bf{0})}\frac{d\bdz}{|(\bdx-h\bdm)-\bdz|^{n-1}}\, \leq c \frac{(h\sqrt{\cD})^n}{|\bdx-h \bdm|^{n-1}}.
\]
If  $|\bdx-h \bdm| < 2 r_0 h \sqrt{\cD}$, by using \eqref{maindis} we obtain that
\[
\int_{B_{r_0 h \sqrt{\cD}}({\bf 0})}
\frac{d\bdz}{|(\bdx-h\bdm)-\bdz|^{n-1}}\leq c (h \sqrt{\cD}).
\]
We deduce that, for all $\bdx \in\R^n$,
\[
\int_{B_{r_0 h \sqrt{\cD}}({\bf 0})}
\frac{1}{|(\bdx-h\bdm)-\bdz|^{n-1}}d\bdz
\leq 
c \frac{(h\sqrt{\cD})^n}{(|\bdx-h \bdm|+h\sqrt{\cD})^{n-1}}.
\]

Therefore
\[
\sum_{h\bdm\in S_h}| \widetilde{f}(h\bdm)|
\int_{{T}\cap B_{r_0 h \sqrt{\cD}}}\!\!\!\left|\kappa_{\l}((\bdx-h \bdm)-\bdz)\!-\!\kappa_{\l}((\bdx-h \bdm)-\bdy(\bdz)\right|
|\eta\left(\!\frac{\bdy(\bdz)}{h\sqrt{\cD}}\right)|d\bdz\]
\[
\leq c_1\, (h \sqrt{\cD})^2 \sum_{\substack{h\bdm\in S_h}}  |\widetilde{f}(h\bdm)| \int_{ B_{r_0 h \sqrt{\cD}}(\bf{0})}\frac{1}{|(\bdx-h\bdm)-\bdz|^{n-1}}\,d\bdz\]
  \[
 \leq C\, ||{f}||_{\infty}(h\sqrt{\cD})^{n+2} \sum_{\substack{h\bdm\in S_h }}
  \frac{1}{(|\bdx-h \bdm|+h \sqrt{\cD})^{n-1}}
  \]
 \[
\leq  C\, ||{f}||_{\infty} (h \sqrt{\cD})^2 \int_{S_h  }
  \frac{d \bdy }{(|\bdx-\bdy|+h\sqrt{\cD})^{n-1}}
\]
\[
\leq  C\, ||{f}||_{\infty}  (h \sqrt{\cD})^3 |\log (\max (h \sqrt{\cD},dist(\bdx,\partial \Omega))|.
\]
The last inequality and \eqref{estimateh2} complete the proof.

\end{proof}

\section{Implementation and numerical examples}
\setcounter{equation}{0}

 From Lemma \ref{halfspace} and Remark \ref{remark2} we obtain the 
following one dimensional
integral representation for the integrals in  \eqref{cub2}
if $h \bdm \in \Omega$ and $d(h \bdm, \partial \Omega) \ge rh \sqrt{\cD}$
\begin{align*}
&\cD^{-n/2}\int_{\R^n}	\kappa_{\l}(\bdx-\bdy)\,
\eta_{2M}\left(\!\frac{\bdy-h\bdm}{h\sqrt{\cD}}\right)d\bdy\\
&=\frac{\l^{n/2-1}}{(2\pi\cD)^{n/2}} \int_{\R^n}
\frac{K_{n/2-1}(\l|\bdx-\bdy|)}{|\bdx-\bdy|^{n/2-1}}\,
\eta_{2M}\left(\!\frac{\bdy-h\bdm}{h\sqrt{\cD}}\right)d\bdy\\
&=\frac{h^{2}}{\cD^{n/2-1}}\frac{1}{(2\pi)^{n/2}} \int_{\R^n}
\frac{(\l\,h\,\sqrt{\cD})^{n/2-1}}{|r_\bdm-\bdy|^{n/2-1}}{K_{n/2-1}(\
\l h\sqrt{\cD}|r_\bdm-\bdy|)}\,
\eta_{2M}\left(\bdy\right)d\bdy\\
&=
\frac{  h^2\cD}{4(\pi\cD)^{n/2}}\int_0^\infty  {\e^{-  \l^2 \,h^2\,\cD
t/4}\e^{-|r_\bdm|^2/(1+  t)}}  \cP_M(r_\bdm,  t)\,dt
\end{align*}
with
\[
r_\bdm=\frac{\bdx-h\bdm}{h\sqrt{\cD}}.
\]

To apply Lemma \ref{halfspace} for nodes $h \bdm$
in the $r h \sqrt{\cD}$ neighborhood
of $\partial \Omega$
we perform a coordinate
transformation such that the center $h\bdm$ is the
origin in the new coordinates $\bdxi=(\xi_1,...,\xi_n)$ and $T_{h\bdm}=\{\bdxi:
\xi_n> \rho_{h\bdm}\}$, where $\rho_{h\bdm}$
is  the distance of the center $h\bdm$ to $\partial\Omega$ .
Since $\eta_{2M}$ is radial we obtain
after the change of variable $\bdy-h\,\bdm=\omega \bdxi h \sqrt{\cD}$, where
$\omega$ is  the rotation matrix,
\begin{align*}
&\frac{\l^{n/2-1}}{(2\pi\,\cD)^{n/2}} \int_{T_{h\bdm}}
\frac{K_{n/2-1}(\l|\bdx-\bdy|)}{|\bdx-\bdy|^{n/2-1}}\,
\eta_{2M}\left(\!\frac{\bdy-h\bdm}{h\sqrt{\cD}}\right)d \bdy\\
&=\frac{\l^{n/2-1}}{(2\pi\,\cD)^{n/2}}
(h\,\sqrt{\cD})^n
\int_{\zeta_n>{\frac{\rho_{h\bdm}}{h\sqrt{\cD}}}}\frac{K_{n/2-1}(\l
h\,\sqrt{\cD} |\omega^{-1}r_{\bdm}-\bdxi|)}{
|\omega^{-1}r_{\bdm}-\bdxi|^{n/2-1}}\eta_{2M}(\bdxi)d\bdxi\\
&=
\frac{  h^n}{8\pi^{n/2}} \int_0^\infty \hskip-3pt \e^{-  \l^2 h^2\cD
\,t/4}\e^{-|r_{\bdm}|^2/(1+  t)} \\
& \qquad\times
\Big(\erfc\big(F_\bdm(t)\big)
\cP_M(\omega^{-1}r_\bdm,  t)
+\frac{\e^{-F^2_\bdm(t)}}
{\sqrt{\pi} }
\cQ_M\big(\omega^{-1}r_\bdm, 
t,\frac{\rho_{h\bdm}}{h\sqrt{\cD}}\big)\Big)\,dt ,
\end{align*}
where we set
\[
F_\bdm(t) =F(  t,(\omega^{-1}r_{\bdm})_n, 
\frac{\rho_{h\bdm}}{h\sqrt{\cD}}) \, .
\]
Then the cubature formula \eqref{cub2} reduces to the computation of 
integrals over the
half-line
\begin{align*}
\widetilde{\cK}_{h,\l}\widetilde{f}(\bdx)&=\frac{   h^2\cD}{4(\pi\cD)^{n/2}}
\sum_{\substack{d(h\bdm,\partial \Omega)\geq
r\,h\sqrt{\cD}\\h\bdm \in \Omega}}{f}(h\bdm)  \int_0^\infty
{\e^{-  \l^2 h^2\,\cD t/4}\e^{-|r_\bdm|^2/(1+  t)}}  \cP_M(r_\bdm,  t)\,dt
\\&+\frac{  h^n}{8\pi^{n/2}}\sum_{\substack{d(h\bdm,\partial
\Omega)<r\,h\sqrt{\cD}}}\widetilde{f}(h\bdm)
\int_0^\infty \e^{-  \l^2 h^2\cD \,t/4}\e^{-|r_{\bdm}|^2/(1+  t)}\\
& \qquad\times \Big(\erfc\big(F_\bdm(t)\big)
\cP_M(\omega^{-1}r_\bdm,  t)
+\frac{\e^{-F^2_\bdm(t)}}
{\sqrt{\pi} }
\cQ_M(\omega^{-1}r_\bdm,  t,\frac{\rho_{h\bdm}}{h\sqrt{\cD}})\Big)\,dt.
\end{align*}

For computing at the grid $\{ h \bdk\}$ we introduce
$\eta_{\bdk-\bdm}=(\omega^{-1}(\bdk-\bdm))_n$ as the projection of
$\bdk-\bdm$ onto the  normal to the tangential plane, i.e. the $\xi_n-$axis,
and set
\begin{align*}
&a_{\bdk}^{(M)}= \frac{1}{4 }\int_0^\infty  \e^{-  \l^2h^2 \, \cD t/4}
\e^{-|\bdk|^2/(\cD(1+  t))}
\cP_M(\frac{\bdk}{\sqrt{\cD}},  t)\,dt\\
&b_{\bdk,\bdm}^{(M)}=\frac{1}{8}\int_0^\infty \e^{-  \l^2 h^2 \, \cD t/4}
\e^{-|\bdk-\bdm|^2/(\cD(1+  t))} \\
& \qquad\times \Big(\erfc\big(F_{\bdk,\bdm}(t)\big)
\cP_M(\frac{\bdk-\bdm}{\sqrt{\cD}},  t)
+\frac{\e^{-F^2_{\bdk,\bdm}(t)}}
{\sqrt{\pi} }
\cQ_M(\omega^{-1}\frac{\bdk-\bdm}{\sqrt{\cD}}, 
t,\frac{\rho_{h\bdm}}{h\sqrt{\cD}})\Big)\,dt,
\end{align*}
where
\[
F_{\bdk,\bdm}(t)= \sqrt{\frac{1+  t}{  t}} \Big( 
\frac{\rho_{h\bdm}}{h\sqrt{\cD}} -
\frac{\eta_{\bdk-\bdm}}{\sqrt{\cD}(1+  t)} \Big).
\]
Hence a  cubature for the volume potential  at the uniform grid
$\{h\bdk\}$ is
\begin{equation}\label{cub4}
\cK_{h,\l}^{(M)} f(h\bdk)=\frac{h^2}{\pi^{n/2}\,\cD^{n/2-1}}\hspace{-5pt}
\sum_{\substack{d(h\bdm,\partial \Omega)\geq
r\,h\sqrt{\cD}\\h\bdm \in
\Omega}}\hspace{-5pt}{f}(h\bdm)\,a_{\bdk-\bdm}^{(M)} + \frac{
h^n}{\pi^{n/2}}\hspace{-5pt}\sum_{\substack{d(h\bdm,\partial
\Omega)<r\,h\sqrt{\cD}}}\hspace{-5pt}\widetilde{f}(h\bdm)b_{\bdk,\bdm}^{(M
)}.
\end{equation}

We transform the one-dimensional integral representation of $a_{\bdk-\bdm}^{(M)}$ and $b_{\bdk,\bdm}^{(M
)}$ to integrals over $\R$  with
  doubly exponentially decaying integrands by making the substitutions
\begin{equation}\label{wald}
t=\e^\xi,\quad \xi=\alpha (\sigma+\e^\sigma),\quad \sigma=\beta (u-\e^{-u}),
\end{equation}
with certain positive constants $a$ and $b$ and the computation is based on the classical trapezoidal rule.
We get after the substitutions
\begin{align*}
&a_{\bdk}^{(M)}= \frac{1}{4 }\int_{-\infty}^\infty  \e^{-  \l^2h^2 \, \cD \Phi(u)/4}
\e^{-|\bdk|^2/(\cD(1+  \Phi(u)))}
\cP_M(\frac{\bdk}{\sqrt{\cD}},  \Phi(u))\,\Phi'(u)du\\
&b_{\bdk,\bdm}^{(M)}=\frac{1}{8}\int_{-\infty}^\infty \e^{-  \l^2 h^2 \, \cD \Phi(u)/4}
\e^{-|\bdk-\bdm|^2/(\cD(1+  \Phi(u)))} \\
& \qquad\times \Big(\erfc\big(F_{\bdk,\bdm}(\Phi(u))\big)
\cP_M(\frac{\bdk-\bdm}{\sqrt{\cD}}, \Phi(u))
+\frac{\e^{-F^2_{\bdk,\bdm}(\Phi(u))}}
{\sqrt{\pi} }
\cQ_M(\omega^{-1}\frac{\bdk-\bdm}{\sqrt{\cD}}, 
\Phi(u),\frac{\rho_{h\bdm}}{h\sqrt{\cD}})\Big)\,\Phi'(u)du,
\end{align*}
where we set
\begin{align*}
\Phi(u)&=\exp(\alpha \beta (u-\exp(-u)) +\alpha \exp(\beta(u-\exp(-u)))),
\\
\Phi'(u)&=\Phi(u) \alpha \beta (1+\e^{-u})(1+\exp(\beta(u-\exp(-u)))).
\end{align*}

The quadrature with the  trapezoidal
rule with step size $\tau$
 gives
\begin{align*}
&a_{\bdk}^{(M)}\approx \frac{\tau}{4 } \sum_{s=-N_0}^{N_1}
\e^{-\l^2\,h^2 \, \cD\,
\Phi(s\,\tau)/4}\e^{-|\bdk|^2/(\cD(1+\Phi(s\,\tau)))}
\cP_M(\frac{\bdk}{\sqrt{\cD}},\Phi(s\,\tau))\,\Phi'(s\,\tau)\, ,\\
&b_{\bdk,\bdm}^{(M)}\approx\frac{\tau}{8} \sum_{s=-N_0}^{N_1}
\e^{-\l^2\,h^2\cD
\,\Phi(s\,\tau)/4}\e^{-|\bdk-\bdm|^2/(\cD(1+\Phi(s\,\tau)))}\Phi'(s\,\tau) \times\\
&
\Big(\erfc\big(F_{\bdk,\bdm}(\Phi(s\,\tau))\big)
\cP_M\big(\frac{\bdk-\bdm}{\sqrt{\cD}},\Phi(s\,\tau)\big)
+\frac{\e^{-F^2_{\bdk,\bdm}(\Phi(s\,\tau))}}{\sqrt{\pi} }
\cQ_M\big(\omega^{-1}\frac{\bdk-\bdm}{\sqrt{\cD}},\Phi(s\,\tau),
\frac{\rho_{h\bdm}}{h\sqrt{\cD}}\big)\Big).
\end{align*}

We have tested	formula \eqref{cub4} in the ellipse
$$\Omega=\{\bdx\in\R^2:\frac{x_1^2}{a^2}+\frac{x_2^2}{b^2}\leq 1\},\quad
0<b\leq a.$$
If $u(\bdx)=0$ and $\nabla u(\bdx)=0$ on $\partial \Omega$ and $0$ outside
the ellipse, then the density
\[
f(\bdx)=(-\Delta+\l^2)u(\bdx),\quad \bdx \in \Omega
\]
satisfies  
\[
\cK_\l f(\bdx)
=\left\{
\begin{tabular}{cc}
$u(\bdx)$& $\bdx\in\Omega$\\
0& $\bdx\notin \Omega$
\end{tabular}
\right.
\]
We have tested the approximation of the potential \eqref{vol}  with the density
\begin{equation}\label{f}
f(\bdx)=(-\Delta
+\l^2)\sin\left(1-\frac{x_1^2}{a^2}-\frac{x_2^2}{b^2}\right)^2
\end{equation}
which provides the exact value
\[
\cK_\l f(\bdx)=
\left\{
\begin{tabular}{cc}
$\sin\left(1-\frac{x_1^2}{a^2}-\frac{x_2^2}{b^2}\right)^2$&
$\bdx\in\Omega$\\
0& $\bdx\notin \Omega$
\end{tabular}
\right.
\]
and we assumed $f$ itself as the extension outside the ellipse.

The results in Table \ref{table1} show the accuracy of the method. They
were obtained with the parameters
$\alpha=4$, $\beta=2$ in \eqref{wald} and the quadrature parameters
$\tau=0.01$,
$N_0=-80$, $N_1=100$. We assumed $h=2^{-7}$  and different $a,b$ and $\l$.
We chose $r=6$ and $\cD=3$, which gives the saturation error less than
$10^{-10}$. 

In Table \ref{table7d4}  and  \ref{table8d4}  we report on the relative
errors and approximation rates for  $\cK_\l f(0.5,0)$ and  $\cK_\l
f(0.25,0.25)$, respectively, for  different $a$ and $b$, and $\l=1$.
  The approximated values are computed  by the cubature formulas
$\cK_{\l,h}^{(M)}$ in \eqref{cub4} for $M=1,2,3$, with $\cD=4$.	We have
chosen
$\alpha=4$, $\beta=2$ in \eqref{wald} and the quadrature parameters
$\tau=0.006$,
$N_0=-160$, $N_1=200$. We obtain rate of convergence $\cO(h^{2M})$  although
Theorem \ref{mainTheorem}. This result comes from  the structure of the error $\cE=|\widetilde{\cK}_{\l,h}\widetilde{f}(\bdx)-\cK_{\l}f(\bdx)|$ which is given by  $\cE=\cE_1+\cE_2$ where 
$\cE_1
=\cO(\epsilon +h^{2M})$ and 
 $\cE_2
 =\cO(h^2)$. Numerical experiments show that  $\cE_2<<\cE_1$.

Analogous  results  are obtained in Table \ref{table8}
for  the approximation of $\cK_\l g(0,0)=1$, with
\begin{equation}\label{g}
g(x_1,x_2)=(-\Delta+\l^2)(1-\frac{x_1^2}{a^2}-\frac{x_2^2}{b^2})^2/ (1+|\bdx|^2)^{-1}.
\end{equation}
To check the rate of convergence we considered the function
\begin{equation}\label{oscill}
f(x_1,x_2)=(-\Delta +\l^2) \cos (30\,\pi\,x_1)\cos(30\,\pi\, x_2).
\end{equation}
The results in Table \ref{table11} show that ${\cK}_{\l,h}^{(M)} \widetilde{f}(\bdx)$ approximates $\cK_\l f (\bdx)$ with the predicted approximation rate $2$, for $M=1,2,3$, in agreement with Theorem \ref{mainTheorem}.

\begin{table}[p]
\begin{scriptsize}
\begin{center}
$a=b=1.5$\\[1mm]
\begin{tabular}{cc|c||cc||cc} \hline
       &       & 	   &  $\l^2=0.2$		&
        &  $\l^2=2$ 	      &\\ \hline
  $x_1$ &$x_2$& exact & approximation &  error  & approximation & 
error\\ \hline
  0.00 &  0.00  &0.8414709848 &0.8414709850 &  0.258E-09    &0.8414709848 &
  0.470E-10  \\
   0.25	&  0.00  &0.8106234643 &0.8106234645 &	0.267E-09    &0.8106234643
&   0.481E-10  \\
   0.50	&  0.00  &0.7104401615 &0.7104401617 &	 0.300E-09    &0.7104401615
&  0.519E-10 \\
   0.75	&  0.00  &0.5333026735 &0.5333026737 &	 0.396E-09  &0.5333026736 &
   0.687E-10  \\
   1.00	&  0.00  &0.3037650630 &0.3037650633 &	0.718E-09  &0.3037650631 &
0.156E-09 \\
   1.25	&  0.00  &0.0932286160 &0.0932286162 &	 0.248E-08    &0.0932286161
&   0.711E-09 \\
0.25 &	0.25&0.7783135137 &0.7783135139 &  0.277E-09  &0.7783135138 &
0.492E-10  \\
   0.50	&  0.50  &0.5687113034 &0.5687113037 &	 0.371E-09  &0.5687113035 &
  0.635E-10  \\
   0.75	&  0.75  &0.2474039593 &0.2474039595 & 0.895E-09  &0.2474039593 &
0.210E-09  \\
   1.00	&  1.00  &0.0123453654 &0.0123453656 &	 0.187E-07  &0.0123453655 &
   0.581E-08 \\
\hline
\end{tabular}\\[1mm]
$a=1.5,b=1$\\[1mm]
\begin{tabular}{cc|c||cc||cc} \hline
      &	      & 	   &  $\l^2=0.2$		&
	&  $\l^2=2$	       &\\ \hline
  $x_1$ &$x_2$& exact & approximation &  error &  approximation & 
error\\ \hline
  0.00 &  0.00  &0.8414709848 &0.8414709867 &   0.219E-08   &0.8414709852 &
  0.519E-09  \\
   0.25	&  0.00  &0.8106234643 &0.8106234661 &	 0.225E-08  &0.8106234647 &
  0.528E-09 \\
   0.50	&  0.00  &0.7104401615 &0.7104401632 &	 0.248E-08   &0.7104401619
&  0.566E-09 \\
   0.75	&  0.00  &0.5333026735 &0.5333026752 &	 0.312E-08  &0.5333026739 &
  0.679E-09 \\
   1.00	&  0.00  &0.3037650630 &0.3037650646 &	 0.511E-08   &0.3037650634
&  0.108E-08 \\
   1.25	&  0.00  &0.0932286160 &0.0932286174 &	 0.154E-07  &0.0932286163 &
   0.321E-08 \\
  0.25  &  0.25	&0.7363058387 &0.7363058405 &	0.247E-08 &0.7363058391 &
0.580E-09 \\
   0.50	&  0.50  &0.3969386023 &0.3969386041 &	 0.457E-08  &0.3969386028 &
   0.117E-08  \\
   0.75	&  0.75  &0.0351490085 &0.0351490104 &	 0.541E-07   &0.0351490091
& 0.182E-07  \\
\hline
\end{tabular}
\\[1mm]
$a=1.5,b=0.5$\\[1mm]
\begin{tabular}{cc|c||cc||ccc} \hline
      &	      & 	   &  $\l^2=0.2$		&
       &  $\l^2=2$		     &\\ \hline
  $x_1$ &$x_2$& exact & approximation &  error &  approximation & 
error\\ \hline
  0.00  &  0.00	&0.8414709848 &0.8414712254 &	0.286E-06  &0.8414710692 &
0.100E-06  \\
   0.25	&  0.00  &0.8106234643 &0.8106237000 &	 0.291E-06  &0.8106235459 &
  0.101E-06  \\
   0.50	&  0.00  &0.7104401615 &0.7104403833 &	 0.312E-06  &0.7104402355 &
   0.104E-06  \\
   0.75	&  0.00  &0.5333026735 &0.5333028743 &	 0.376E-06  &0.5333027358 &
   0.117E-06  \\
   1.00	&  0.00  &0.3037650630 &0.3037652384 &	 0.577E-06   &0.3037651117
&   0.160E-06 \\
   1.25	&  0.00  &0.0932286160 &0.0932287650 &	 0.160E-05   &0.0932286511
&  0.377E-06  \\
  0.25  &  0.25	&0.4982722935 &0.4982725309 &	0.476E-06  &0.4982723784 &
  0.170E-06  \\
  \hline
\end{tabular}
\end{center}
\caption{
Exact and approximated values of  $\cK_\l f(x_1,x_2)$ and relative error
using $\cK_{\l,h}^{(3)} f(x_1,x_2)$,
when $h=2^{-7}$ and $f$ is given in \eqref{f} }\label{table1}
\end{scriptsize}
\end{table}

\begin{table}[p]
\begin{scriptsize}
\begin{center}
\begin{tabular}{r|c|c|c} \hline
\end{tabular}\\[1mm]
$M=1$\\[1mm]
\begin{tabular}{r|cc|cc|cc} \hline
& $a=1.5$&$b=1.5$ &$a=1.5$&$b=1$  &$a=1.5$& $b=0.5$\\ \hline
$h^{-1} \hspace{-2mm}$ &  error  & rate &  error  & rate & 
error &rate   \\[1pt] \hline
$2^4$  &	0.439E-01 &	&0.968E-01 &	 &0.572E+00 &	     \\
$2^5$&	0.110E-01 &    1.997&0.243E-01 &    1.993&0.167E+00 &	 1.781 \\
$2^6$  & 0.275E-02 &    2.000&0.608E-02 &    2.000&0.419E-01 &	1.990  \\
$2^7$  &0.688E-03 &    2.000&0.152E-02 &	  2.000&0.105E-01 &    2.000	\\
$2^8$ &	0.172E-03 &    2.000&0.380E-03 &    2.000&0.262E-02 &	 2.000 \\
$2^9$ &  0.430E-04 &    2.000&0.950E-04 &    2.000&0.655E-03 &	2.000  \\
\hline
\end{tabular}\\[1mm]
$M=2$\\[1mm]
\begin{tabular}{r|cc|cc|cc} \hline
& $a=1.5$&$b=1.5$ &$a=1.5$&$b=1$  &$a=1.5$& $b=0.5$\\ \hline
$h^{-1} \hspace{-2mm}$ &  error  & rate &  error  & rate & 
error &rate   \\[1pt] \hline
$2^4$  &	0.174E-03 &    &0.114E-02 &	&0.186E+00 &	  \\
$2^5$ &	 0.443E-05 &	5.294&0.500E-05 &    7.831&0.288E-02 &	  6.013 \\
$2^6$& 0.183E-06 &    4.601&0.626E-06 &    2.996&0.218E-04 &	7.044  \\
$2^7$  &0.996E-08 &    4.196 &0.534E-07 &    3.552&0.919E-06 &	4.570	 \\
$2^8$ &	0.600E-09 &    4.053 &0.356E-08 &    3.908&0.922E-07 &	  3.317 \\
$2^9$ &  0.371E-10 &    4.013 &0.226E-09 &    3.979&0.630E-08 &	 3.871	\\
\hline
\end{tabular}\\[1mm]
$M=3$\\[1mm]
\begin{tabular}{r|cc|cc|cc} \hline
& $a=1.5$&$b=1.5$ &$a=1.5$&$b=1$  &$a=1.5$& $b=0.5$\\ \hline
$h^{-1} \hspace{-2mm}$ &  error  & rate &  error  & rate & 
error &rate   \\[1pt] \hline
$2^4$  &	0.719E-04 &   &0.747E-03 &    &0.469E-01 &    \\
$2^5$ &	0.102E-05 &    6.138&0.102E-04 &    6.189&0.177E-02 &	 4.732 \\
 $2^6$ & 0.155E-07 &    6.039&0.153E-06 &    6.066&0.248E-04 &	6.152  \\
$2^7$  &0.241E-09 &    6.010&0.236E-08 &	  6.018&0.373E-06 &    6.057	\\
$2^8$ &	0.376E-11 &    6.001&0.368E-10 &    6.003 &0.577E-08 &	  6.015 \\
$2^9$ &  0.936E-13 &    5.328& 0.507E-12 &    6.180&0.899E-10 &	 6.003	\\
\hline
\end{tabular}\\[1mm]
\caption{ Relative errors and approximation rates
for $\cK_\l f(0.5,0)=0.7104401614873481 $ using $\cK^{(M)}_{\l,h} f(0.5,0)$
in $\Omega=\{\bdx\in\R^2:\frac{x_1^2}{a^2}+\frac{x_2^2}{b^2}\leq 1\}$, $\l=1$ and $f$ in \eqref{f} }\label{table7d4}
\end{center}
\end{scriptsize}
\begin{scriptsize}
\begin{center}
$M=1$\\[1mm]
\begin{tabular}{r|cc|cc|cc} \hline
& $a=1.5$&$b=1.5$ &$a=1.5$&$b=1$  &$a=1.5$& $b=0.5$\\ \hline
$h^{-1} \hspace{-2mm}$ &  error  & rate &  error  & rate & 
error &rate   \\[1pt] \hline
$2^4$  &	0.387E-01 &    &0.955E-01 &	& 0.822E+00 &	   \\
$2^5$ &	0.967E-02 &    2.000& 0.240E-01 &    1.991&0.246E+00 &	  1.738 \\
 $2^6$  & 0.242E-02 &    2.000&0.601E-02 &    1.999&0.623E-01 &	1.983  \\
 $2^7$  &0.604E-03 &    2.000&0.150E-02 &	  2.000& 0.156E-01 &	1.999	\\
$2^8$&	0.151E-03 &    2.000&0.376E-03 &    2.000&0.390E-02 &	 2.000 \\
$2^9$ &  0.378E-04 &    2.000& 0.939E-04 &    2.000&0.974E-03 &	 2.000	 \\
\hline
\end{tabular}\\[1mm]
$M=2$\\[1mm]
\begin{tabular}{r|cc|cc|cc} \hline
& $a=1.5$&$b=1.5$ &$a=1.5$&$b=1$  &$a=1.5$& $b=0.5$\\ \hline
$h^{-1} \hspace{-2mm}$ &  error  & rate &  error  & rate & 
error &rate   \\[1pt] \hline
$2^4$  &	 0.593E-04 &	 &0.139E-02 &	 &0.312E+00 &  \\
$2^5$ &	0.225E-05 &    4.721&0.192E-04 &    6.174&0.586E-02 &	 5.735 \\
$2^6$  & 0.228E-06 &    3.301&0.244E-06 &    6.297&0.114E-03 &	5.682  \\
$2^7$  & 0.156E-07 &    3.870&0.765E-09 &    8.318&0.350E-05 &	5.026	 \\
$2^8$&	0.997E-09 &    3.969&0.177E-09 &    2.113&0.164E-06 &	 4.420 \\
$2^9$ &  0.627E-10 &    3.992&0.145E-10 &    3.612&0.937E-08 &	4.126  \\
\hline
\end{tabular}\\[1mm]
$M=3$\\[1mm]
\begin{tabular}{r|cc|cc|cc} \hline
& $a=1.5$&$b=1.5$ &$a=1.5$&$b=1$  &$a=1.5$& $b=0.5$\\ \hline
$h^{-1} \hspace{-2mm}$ &  error  & rate &  error  & rate & 
error &rate   \\[1pt] \hline
$2^4$  &	0.663E-04 &    &0.762E-03 &	&0.677E-01 &	   \\
$2^5$ &	0.947E-06 &    6.130&0.104E-04 &    6.191&0.281E-02 &	 4.589 \\
 $2^6$  &	0.144E-07 &    6.037&0.156E-06 &    6.067&0.395E-04 &	 6.154	\\
$2^7$  &0.224E-09 &    6.009& 0.240E-08 &    6.018&0.593E-06 &	6.058	 \\
$2^8$&	0.348E-11 &    6.008&0.374E-10 &    6.004&0.917E-08 &	 6.016 \\
$2^9$&  0.117E-12 &    4.896&0.558E-12 &    6.067&0.143E-09 &	6.003  \\
\hline
\end{tabular}\\[1mm]
\caption{ Relative  errors and approximation rates
for $\cK_\l f(0.25,0.25)$ using $\cK^{(M)}_{\l,h} f(0.25,0.25)$ in
$\Omega=\{\bdx\in\R^2:\frac{x_1^2}{a^2}+\frac{x_2^2}{b^2}\leq 1\}$, $\l=1$ and $f$ in \eqref{f}. }\label{table8d4}
\end{center}
\end{scriptsize}
\end{table}
\begin{table}[p]
\begin{scriptsize}
\begin{center}
$M=1$\\[1mm]
\begin{tabular}{r|cc|cc|cc} \hline
& $a=1.5$&$b=1.5$ &$a=1.5$&$b=1$  &$a=1.5$& $b=0.5$\\ \hline
$h^{-1} \hspace{-2mm}$ &  error  & rate &  error  & rate & 
error &rate   \\[1pt] \hline
 $2^2$  &0.415E+00 &&0.648E+00 &&0.368E+01 &     \\
 $2^3$ &0.139E+00 &    1.574&0.216E+00 &	 1.585&0.131E+01 &    1.489  \\
$2^4$ &0.386E-01 &    1.854&0.594E-01 &	 1.862&0.366E+00 &    1.841  \\
$2^5$ &0.992E-02 &    1.959&0.153E-01 &	1.961&0.942E-01 &    1.959 \\
 $2^6$  &0.250E-02 &    1.989&0.384E-02 &	  1.990&0.237E-01 &    1.989  \\
$2^7$ &  0.626E-03 &    1.997&0.962E-03 &    1.997&0.594E-02 &	1.997 \\
$2^8$ &0.157E-03 &    1.999&0.241E-03 &	 1.999&0.149E-02 &    1.999\\
$2^9$ & 0.391E-04 &    2.000&0.602E-04 &	  2.000&0.372E-03 &    2.000  \\
\hline
\end{tabular}\\[1mm]
$M=2$\\[1mm]
\begin{tabular}{r|cc|cc|cc} \hline
& $a=1.5$&$b=1.5$ &$a=1.5$&$b=1$  &$a=1.5$& $b=0.5$\\ \hline
$h^{-1} \hspace{-2mm}$ &  error  & rate &  error  & rate & 
error &rate   \\[1pt] \hline
 $2^2$   &0.134E+00 & &0.210E+00 & &0.141E+01 &   \\
$2^3$ &0.173E-01 &    2.952& 0.254E-01 &	  3.050&0.187E+00 &    2.920   \\
$2^4$  &0.143E-02 &    3.599&0.207E-02 &	 3.616&0.138E-01 &    3.756  \\
$2^5$ &0.968E-04 &    3.882&0.140E-03 &	3.885&0.914E-03 &    3.918 \\
 $2^6$  &0.618E-05 &    3.969&0.895E-05 &	  3.970&0.580E-04 &    3.978   \\
$2^7$  &0.388E-06 &    3.992&0.562E-06 &	  3.992&0.364E-05 &    3.995  \\
$2^8$ &0.243E-07 &    3.998& 0.352E-07 &	  3.998&0.228E-06 &    3.999\\
$2^9$ &0.152E-08 &    3.999&0.220E-08 &	 3.999& 0.142E-07 &    4.000 \\
\hline
\end{tabular}\\[1mm]
$M=3$\\[1mm]
\begin{tabular}{r|cc|cc|cc} \hline
& $a=1.5$&$b=1.5$ &$a=1.5$&$b=1$  &$a=1.5$& $b=0.5$\\ \hline
$h^{-1} \hspace{-2mm}$ &  error  & rate &  error  & rate & 
error &rate   \\[1pt] \hline
$2^2$     &0.495E-01 &&0.703E-01 &&0.577E+00 &    \\
$2^3$ &0.284E-02 &    4.120 &0.394E-02 &	  4.157&0.229E-01 &    4.659  \\
$2^4$  &0.751E-04 &    5.244&0.106E-03 &	 5.221&0.499E-03 &    5.517  \\
$2^5$ &0.138E-05 &    5.768&0.195E-05 &	5.762 &0.903E-05 &    5.789 \\
$2^6$  &0.225E-07 &    5.938&0.318E-07 &	  5.937&0.145E-06 &    5.960  \\
$2^7$  &0.355E-09 &    5.984& 0.502E-09 &    5.984 &0.228E-08 &	 5.990 \\
$2^8$&0.549E-11 &    6.014&0.790E-11 &	 5.992 &0.357E-10 &    5.998\\
$2^9$ & 0.150E-12 &    5.190&0.858E-13 &	  6.524&0.576E-12 &    5.954  \\
\hline
\end{tabular}\\[1mm]
\caption{ Relative errors and approximation rates
for $\cK_\l g(0,0)=1$ using $\cK^{(M)}_{\l,h} g(0,0)$ in
$\Omega=\{\bdx\in\R^2:\frac{x_1^2}{a^2}+\frac{x_2^2}{b^2}\leq 1\}$,  $\l=1$ and $g$ in \eqref{g}. }\label{table8}
\end{center}
\end{scriptsize}
\end{table}
\begin{table}[p]
\begin{scriptsize}
\begin{center}
$M=1$\\[1mm]
\begin{tabular}{r|r|cc|cc|cc} \hline
& &$a=1.5$&$b=1.5$ &$a=1.5$&$b=1$  &$a=1.5$& $b=0.5$\\ \hline
$h^{-1} \hspace{-2mm}$&$x_1$ &	 error  & rate &   error	& rate &   error
&rate	\\[1pt] \hline
$2^8$ &  0.250	&   0.171E-05 && 0.171E-05&&  0.171E-05 &   \\
'' &  0.750  &	 0.395E-05 &&  0.395E-05 &&  0.395E-05 &    \\
'' &  1.250  &	0.268E-05 &&  0.268E-05 &&  0.268E-05 &   \\
'' &  1.500  &	  0.569E-08 &&	0.592E-08&&  0.688E-08 &  \\ \hline
$2^9$&	0.250  &    0.762E-06  & 1.163 &  0.762E-06 &	  1.163&  0.762E-06
&   1.163  \\
''&  0.750  &	0.176E-05 & 1.163&  0.176E-05 &     1.163 &  0.176E-05 &
1.162  \\
''&  1.250  &	 0.120E-05 & 1.162&  0.120E-05 &      1.162 &  0.120E-05 &
  1.162 \\
''&  1.500  &	 0.140E-07 &  -1.302&  0.144E-07 &    -1.284 &	0.163E-07 &
   -1.241  \\  \hline
$2^{10}$&  0.250  &   0.217E-06 &  1.811&0.217E-06 &  1.811 &  0.217E-06 &
1.811 \\
''&  0.750  &	0.503E-06 &  1.811&0.503E-06 &	1.811 &  0.503E-06 &
1.811 \\
''&  1.250  &	 0.341E-06 & 1.811&0.341E-06 &	1.811&	0.341E-06 &   1.811
\\
''&  1.500  &  0.460E-08 & 1.609&  0.467E-08 & 1.627&  0.504E-08 &   1.690
\\  \hline
$2^{11}$&  0.250  &  0.561E-07 &  1.954& 0.561E-07 &	 1.954&  0.561E-07
&   1.954 \\
''& 0.750  & 0.130E-06 &  1.954&0.130E-06 &  1.954 &  0.130E-06 & 1.954 \\
''&  1.250  &	0.881E-07 &  1.954&0.881E-07 &	1.954 &  0.881E-07 &
1.954\\
''&  1.500  &  0.122E-08 &  1.915&  0.123E-08 & 1.926&	0.128E-08 & 1.975
\\  \hline
\end{tabular}\\[1mm]
$M=2$\\[1mm]
\begin{tabular}{r|r|cc|cc|cc} \hline
& &$a=1.5$&$b=1.5$ &$a=1.5$&$b=1$  &$a=1.5$& $b=0.5$\\ \hline
$h^{-1} \hspace{-2mm}$&$x_1$ &	 error  & rate &   error	& rate &   error
&rate	\\[1pt] \hline
$2^8$ &  0.250	& 0.316E-05  &&0.316E-05 &  &  0.316E-05 &    \\
'' &  0.750  &0.731E-05  & &0.731E-05 &  &  0.731E-05 &  \\
'' &  1.250  &	 0.497E-05&& 0.497E-05 & &  0.496E-05 &     \\
'' &  1.500  & 0.534E-07   & &	0.555E-07 &  &	0.644E-07 & \\	\hline
$2^9$&	0.250  &0.897E-06 &1.816&  0.897E-06 &	  1.816&  0.897E-06 &
1.816 \\
''&  0.750  &0.208E-05 &1.816&	0.208E-05 &   1.816& 0.208E-05 &
1.816\\
''&  1.250  &0.141E-05 &1.816 &  0.141E-05 &   1.816 &	0.141E-05 &
1.815\\
''&  1.500  &0.194E-07	&    1.464 &  0.196E-07 &  1.503 &  0.210E-07  &
1.619	\\  \hline
$2^{10}$&  0.250  &0.227E-06  &1.986 &	0.227E-06 &	1.986 &  0.227E-06
&     1.986\\
''&  0.750  &0.524E-06	 &1.986&  0.524E-06 &	1.986 &  0.524E-06 &
1.986 \\
''&  1.250  &0.356E-06	 &1.986&  0.356E-06 &	1.986&	0.356E-06 &
1.986 \\
''&  1.500  &0.494E-08	&1.970&  0.495E-08 &1.984&  0.504E-08  &    2.056
\\   \hline
$2^{11}$&  0.250  &0.567E-07&1.999&  0.567E-07 &  1.999 &  0.567E-07 &
1.999 \\
''& 0.750  &0.131E-06 &1.999&  0.131E-06 &    1.999 &  0.131E-06  &
1.999 \\
''&  1.250  &0.890E-07 &1.999&	0.890E-07 &   1.999 &  0.890E-07 &    1.999
   \\
''&  1.500  &0.124E-08 & 1.999&  0.124E-08 &  2.001 &  0.124E-08  &
2.022 \\
\hline
\end{tabular}\\[1mm]
$M=3$\\[1mm]
\begin{tabular}{r|r|cc|cc|cc} \hline
& &$a=1.5$&$b=1.5$ &$a=1.5$&$b=1$  &$a=1.5$& $b=0.5$\\ \hline
$h^{-1} \hspace{-2mm}$&$x_1$ &	 error  & rate &   error	& rate &   error
&rate	\\[1pt] \hline
$2^8$ &  0.250	&0.356E-05& &  0.356E-05 &&  0.355E-05 & \\
'' &  0.750  &0.823E-05&&  0.823E-05 &	&  0.823E-05 &	\\
'' &  1.250  &	0.559E-05&&  0.559E-05 &  &  0.559E-05 &   \\
'' &  1.500  &0.751E-07& &  0.770E-07 & &  0.862E-07 &\\ \hline
$2^9$&	0.250  &  0.906E-06 &	  1.972&  0.906E-06 &	  1.972 &
0.906E-06 &   1.972 \\
''&  0.750  &	0.210E-05 &    1.972&  0.210E-05 &    1.972 &  0.210E-05 &
   1.972\\
''&  1.250  &  0.142E-05 &     1.972 &	0.142E-05 &   1.972 &  0.142E-05 &
   1.972 \\
''&  1.500  &  0.198E-07 &     1.925 &	0.198E-07 &  1.959&  0.203E-07 &
2.087 \\ \hline
$2^{10}$&  0.250  & 0.227E-06 &   1.999 &  0.227E-06 &	  1.999 &
0.227E-06 &	 1.999 \\
''&  0.750  &  0.525E-06 &   1.999&  0.525E-06 &    1.999&  0.525E-06 &
1.999\\
''&  1.250  &  0.356E-06 &    1.999 &  0.356E-06 &    1.999 &  0.356E-06 &
   1.999  \\
''&  1.500  &  0.495E-08 & 1.999&  0.495E-08 &	  2.002 &  0.495E-08 &
2.035\\ \hline
$2^{11}$&  0.250  &  0.567E-07 &      2.000&  0.567E-07 &   2.000&
0.567E-07 &	 2.000	\\
''& 0.750  & 0.131E-06 &      2.000&  0.131E-06 &   2.000&  0.131E-06 &
2.000 \\
''&  1.250  &  0.891E-07 &     2.000&  0.891E-07 &    2.000&  0.891E-07 &
  2.000 \\
''&  1.500  &  0.124E-08 &     2.000&  0.124E-08 &  2.000  &  0.124E-08 &
  2.001\\
\hline
\end{tabular}\\[1mm]
\caption{ Absolute errors and approximation rates
for $\cK_1 f(x_1,0)$  using $\cK^{(M)}_{1,h} f(x_1,0)$ in
$\Omega=\{\bdx\in\R^2:\frac{x_1^2}{a^2}+\frac{x_2^2}{b^2}\leq 1\}$, with $f$ given in
\eqref{oscill}. }\label{table11}
\end{center}
\end{scriptsize}
\end{table}

\end{document}